\documentclass[a4paper,10pt,notitlepage,reqno]{article}
\usepackage[english]{babel}
\usepackage{amssymb,amsthm,amsmath} 
\usepackage{mathtools}
\usepackage{mathrsfs}
\usepackage{enumerate}
\usepackage{authblk}
\usepackage{cite}
\usepackage[symbol]{footmisc}
\usepackage{color}
\usepackage{geometry}
\usepackage[colorlinks=true, linkcolor=red, citecolor=blue, urlcolor=blue,hyperfootnotes=false]{hyperref}
\theoremstyle{plain} 
\newtheorem{theorem}{Theorem}[section]
\newtheorem{lemma}{Lemma}[section]

\theoremstyle{definition}
\newtheorem{definition}{Definition}[section]

\theoremstyle{remark}
\newtheorem{remark}{Remark}[section]

\DeclareMathOperator{\supp}{supp}

\newcommand{\R}{\mathbb{R}}

\newcommand{\N}{\mathbb{N}}
\newcommand{\normeq}[1]{{\left\vert\kern-0.25ex\left\vert\kern-0.25ex\left\vert #1 
    \right\vert\kern-0.25ex\right\vert\kern-0.25ex\right\vert}}

{\left\lbrace\begin{array}{r@{\hspace{1mm}}ll}}%
{\end{array}\right.}

\title{\textbf{On uniqueness of solutions to the surface electromigration equation}}

\author[1]{Lucrezia Cossetti} 
\author[1,2]{Luca Fanelli}		
\author[3]{Felipe Linares}

\affil[1]{Ikerbasque \& Universidad del Pa\'is Vasco/Euskal Herriko Unibertsitatea, UPV/EHU, Aptdo. 644, 48080, Bilbao, Spain; lucrezia.cossetti@ehu.eus; luca.fanelli@ehu.es}

\affil[2]{BCAM-Basque Center for Applied Mathematics, Mazarredo, 14 E48009 Bilbao, Spain}

\affil[3]{IMPA, Estrada Dona Castorina 110, Rio de Janeiro 22460-320, Brazil; linares@impa.br}

\begin{document}

\date{\small 7 March 2023}

\maketitle



\begin{abstract}
	\noindent
	In this paper we investigate on uniqueness properties of solutions to the surface electromigration (SEM) equation, which is a generalisation of the more classical Zakharov-Kuznetsov equation of plasma physics with non-local perturbation terms. We will show that if the difference of two solutions has a sufficiently strong spatial decay at two different instants of time, then the two solutions coincide on the whole interval of time. 
\end{abstract}



\section{Introduction}

In this paper we investigate on uniqueness properties of solutions to the initial value problem for the surface electromigration (SEM) equation
\begin{equation}\label{eq:zkp}
	\begin{cases}
		\partial_t u + \partial_x \Delta u + \frac{1}{2} (u \partial_x u + \partial_x u \partial_x \phi + \partial_y u \partial_y \phi)=0, \quad (x,y)\in \R^2,\, t\in [0,1],\\
		\Delta \phi=\partial_x u.
	\end{cases}
\end{equation}	
Here $u=u(x,y,t)$ and $\phi=\phi(x,y,t)$ are real-valued functions and $\Delta$ denotes the two-dimensional Laplacian operators.

The surface electromigration (SEM) phenomenon has been largely studied in the last decades ~\cite{GuMa1998,GuMa1999,ScKr1997,Bradley1999} and it can be described as follows: let a constant, unidirectional, electric current pass trough a solid metal film (nanoconductor), the collisions between the current-carrying electrons and the metal atoms at the surface of the film lead to drift of these atoms. This effect can cause motion and/or deformation of the solid metal surface, a response which is commonly known as surface electromigration (SEM) phenomenon. In~\cite{Bradley2001} it was shown that, in the limit of high applied currents, \emph{i.e.} considering the surface electromigration effect as the predominant cause of atom motion, the propagation of the mentioned surface disturbance is fruitfully described by system~\eqref{eq:zkp} for regimes of small amplitude and large width waves.  More specifically in~\eqref{eq:zkp} $u$ plays the role of the surface displacement and $\phi$ is used to denote the electrostatic potential on the surface conductor.    

If one sets $\phi$ in~\eqref{eq:zkp} to zero, then the first equation in system~\eqref{eq:zkp} reduces to the well known Zakharov-Kuznetsov (ZK) equation
\begin{equation}\label{eq:zk}
\partial_t u + \partial_x \Delta u + u\partial_x u=0,
\quad (x,y)\in \R^2,\, t\in [0,1].
\end{equation}   
Thus, system~\eqref{eq:zkp} differs from the ZK equation through the coupling with the potential equation for $\phi.$ 

The IVP associated to the SEM equation~\eqref{eq:zkp} was investigated by Linares, Pastor and Scialom in~\cite{LPS2021}. In this paper they established a local well-posedness theory in $H^s(\R^2),$ $s>1/2.$ Their result made use of the available local-well posedness theory for ZK equation~\eqref{eq:zk} that in the recent years has seen a remarkable development. The first result for ZK was due to Faminskii~\cite{Faminskii1995} who proved local well posedness in $H^s(\R^2),$ $s\geq 1.$ Later, in~\cite{LP2009}, Linares and Pastor proved an analogous result lowering the Sobolev exponent to $s\geq 3/4.$ Further improvements were established independently by Grünrock and Herr in~\cite{GH2014} and Molinet and Pilot in~\cite{MP2015}. In these works the author pushed down the Sobolev index to $s> 1/2.$ Recently, Kinoshita~\cite{Kinoshita2021} established the sharp local well-posedness in $H^s(\R^2),$ $s>-1/4$ (see also~\cite{LS2009,MP2015,RV2012,HK2020} for the corresponding results in the higher-dimensional case).

Asymptotic behaviour of solutions to ZK have been studied in~\cite{MMPP2021,FLP2012,CMPS2016,KRS2021,Valet2020,deBouard1996,FHRY2018,FHRY2020}. Propagation of regularity has been considered in~\cite{LP2018,MR2023,Mendez2020}. See also the recent review~\cite{Munoz2022}.
Another aspect that has strongly attracted the interest of the dispersive community is the unique continuation principle. Roughly speaking, proving a unique continuation principle for dispersive equations means providing sufficient conditions on the (spatial) decay of the difference $u_1-u_2$ of two solutions of the model at two different times which guarantee that the solutions coincide in the entire interval of time. In connection with the ZK equation this topic was already investigated by Bustamante, Isaza and Mej\'ia in~\cite{BIM2013} in the 2D case and later by the authors of this paper in~\cite{CFL2019}. Later the strategy developed in~\cite{BIM2013} was generalized by Bustamante, Jim\'enez and Mej\'ia in~\cite{BJM2018} to establish a unique continuation principle for the ZK equation in the 3D case.

In this note we are interested to extending these uniqueness results to the surface electromigration (SEM) equation~\eqref{eq:zkp}. Before stating our main result we rewrite conveniently system~\eqref{eq:zkp} as a single equation: defining the operator $\mathcal{L}:=(\Delta)^{-1}\partial_x,$ then system~\eqref{eq:zkp} is equivalent to
\begin{equation}\label{eq:zknl}
	\partial_t u + \partial_x \Delta u + 
	\frac{1}{2}\Big(u\,\partial_x u + \partial_x u\, \partial_x \mathcal{L}(u) + \partial_y u\, \partial_y \mathcal{L}(u)\Big)=0,
	\quad (x,y)\in \R^2,\, t\in [0,1].
\end{equation}
Our uniqueness result reads as follows.
\begin{theorem}\label{thm:main}
	Suppose that for some small $\varepsilon>0$
	\begin{equation*}
		u_1, u_2\in C\big([0,1];H^4(\R^2)\cap L^2((1+ x^2+y^2)^{4/3 + \varepsilon} dx dy)\big)\cap C^1\big([0,1]; L^2(\R^2)\big),
	\end{equation*}
	are solutions to~\eqref{eq:zknl}. Then there exists a universal constant $a_0>0$ such that if for some $a>a_0$
	\begin{equation}\label{eq:exp_cond}
		u_1(0)-u_2(0), u_1(1)-u_2(1)\in L^2\big(e^{a(x^2+y^2)^{3/4}} dx dy\big),
	\end{equation}
	then $u_1=u_2.$
\end{theorem}
The proof of Theorem~\ref{thm:main} follows the scheme developed in~\cite{EKPV2007} for the KdV equation and it is based upon the comparison of two types of estimates, a lower bound, which comes as a consequence of a Carleman estimate for suitable compactly supported functions and an upper bound for the $H^2$ norm of the solution which follows from the exponential decay assumed in~\eqref{eq:exp_cond} for the data.  
We stress that for both the lower and the upper bound we will need to adapt/generalise the results available for the standard ZK equation, among those a better suited Carleman estimate will be needed. Moreover, additional difficulties arise when dealing with the new non-local terms characterising the SEM model.   

\subsection*{Preliminary facts and notations}
We start this preliminary section writing more explicitly the non-local terms, $\partial_x \mathcal{L}$ and $\partial_y \mathcal{L},$ in~\eqref{eq:zknl}. As one can see, they are essentially  square of the Riesz transform. Indeed, denoting with $(\xi, \eta)$ the Fourier variables corresponding to the physical space variables $(x,y)$, one has
	\begin{equation}\label{eq:non-localF}
		\mathcal{F}(\partial_x \mathcal{L}g)(\xi, \eta)
		=\frac{\xi^2}{\xi^2 + \eta^2} \mathcal{F}(g)(\xi, \eta),
		\qquad 
		\mathcal{F}(\partial_y \mathcal{L}g)(\xi, \eta)
		=\frac{\xi\eta}{\xi^2 + \eta^2} \mathcal{F}(g)(\xi, \eta).
	\end{equation}
	We recall that the Riesz transform $\mathcal{R}=(\mathcal{R}_1, \mathcal{R}_2,\dots, \mathcal{R}_d)$ is defined for any $g\in L^2(\R^d),$ in Fourier space, by
	\begin{equation}\label{eq:Riesz}
		\mathcal{F}(\mathcal{R}_j g)(\zeta)=-i\frac{\zeta_j}{|\zeta|}\mathcal{F}(g)(\zeta),
		\qquad j=1,2,\dots, d.
	\end{equation}
	Thus, from definitions~\eqref{eq:non-localF} and~\eqref{eq:Riesz}, one sees that the operators $\partial_x \mathcal{L}$ and $\partial_y \mathcal{L}$ can be written more compactly in terms of the 2D Riesz transform $\mathcal{R}=(\mathcal{R}_x, \mathcal{R}_y)$ as
	\begin{equation}\label{eq:deLRiesz}
		\partial_x \mathcal{L}=-\mathcal{R}_x^2 ,
		\qquad
		\partial_y \mathcal{L}=-\mathcal{R}_x\mathcal{R}_y.
	\end{equation}
For later purposes we need the definition of the Muckenhoupt class of weights $A_p(\R^d).$
\begin{definition}
	Let $1<p<\infty$ and let $w$ be a measurable non-negative function. We say that $w\in A_p(\R^d)$ if 
	\begin{equation}\label{eq:Ap}
		Q_p(w):=\sup_{B} \left( \frac{1}{|B|} \int_{B} w(x)\, dx\right) \left( \frac{1}{|B|} \int_{B} w(x)^{-\frac{1}{p-1}}\, dx \right)^{p-1}\leq C,
	\end{equation}	 
	where $B$ is any ball in $\R^d$ and $C$ is a constant independent of $B.$
\end{definition}
We recall the following result on boundedness, respectively, weighted boundedness of the Riesz transform (see~\cite{BW1995,Petermichl2008,CF1974}).
\begin{lemma}[Boundedness Riesz transform]
\label{lemma:boundedness}
	Let $1<p<\infty$ and $p'$ such that $1/p+1/p'=1$ and let $w$ be a weight in the $A_p(\R^d)$-class (see definition~\eqref{eq:Ap}). Then, for any $j=1, 2 \dots, d$ the following bounds on the operator norms of the Riesz transform $\mathcal{R}_j$ hold true:
	\begin{align}
		\label{eq:Riesz-boundedness}
		\|\mathcal{R}_j\|_{L^p\to L^p}=\cot\Big( \frac{\pi}{2p^\ast}\Big),
		\qquad &p^\ast:=\max\{p,p'\},\\
		\label{eq:Riesz-weighted-boundedness}
		\|\mathcal{R}_j\|_{L^p(w)\to L^p(w)}\leq c_{p,d} Q_p(w)^r,
		\qquad &r:=\max \{1, p'/p\}
	\end{align}
	Moreover, both the bounds are sharp, \emph{i.e.} the best possible bound is established.
\end{lemma}

\section{Proof of Theorem~\ref{thm:main}}
This section is concerned with the proof of Theorem~\ref{thm:main}. As customarily in proving unique continuation results, this will follow as a consequence of a comparison between suitable lower and upper bounds for the difference $v:=u_1-u_2$ of two solutions $u_1,u_2$ of~\eqref{eq:zknl}. Working with $v$ turns the problem of proving uniqueness for the equation~\eqref{eq:zknl} into the problem of proving the triviality of the solution to the corresponding equation associated to $v,$ namely
\begin{equation}\label{eq:eq_v}
\partial_t v + \partial_x \Delta v
+ \frac{1}{2}
\big(
u_1\partial_x v + \partial_x u_2 v 
+ \partial_x u_1 \partial_x \mathcal{L}(v) + \partial_x \mathcal{L}(u_2)\partial_x v
+ \partial_y u_1 \partial_y \mathcal{L}(v) + \partial_y \mathcal{L}(u_2)\partial_y v
\big)=0.
\end{equation} 
In the following we will consider a more general problem than~\eqref{eq:eq_v}, namely
\begin{equation}\label{eq:eq_v_gen}
	\partial_t v + \partial_x \Delta v + a_1(x,y,t)\partial_x v +  b_1(x,y,t)\partial_y v
	+ a_0(x,y,t)\partial_x \mathcal{L}(v) +  b_0(x,y,t)\partial_y \mathcal{L}(v)
	+ c_0(x,y,t) v=0,
\end{equation}
for suitable variable coefficients $a_1,b_1, a_0,b_0$ and $c_0.$ Clearly~\eqref{eq:eq_v} is a particular case of~\eqref{eq:eq_v_gen}, this can be easily seen setting $a_1=\frac{1}{2}(u_1 + \partial_x \mathcal{L}(u_2)),$ $b_1=\frac{1}{2}\partial_y \mathcal{L}(u_2),$ $a_0=\frac{1}{2}\partial_x u_1,$ $b_0=\frac{1}{2}\partial_y u_1$ and $c_0=\frac{1}{2}\partial_x u_2.$

We start introducing explicitly the quantity that we aim at estimating from below and from above, namely
\begin{equation}\label{eq:A_R(v)}
	A_R(v):= \left( \int_0^1 \int_{Q_R} (|v|^2 + |\partial_x v|^2 + |\partial_y v|^2 + |\partial_x\partial_y v|^2 + |\Delta v|^2) dxdydt \right)^{1/2},
\end{equation}
where $Q_R:=\{(x,y)\colon R-1\leq \sqrt{x^2+y^2}\leq R\}.$

The statements of the lower and upper bounds are contained in the following two theorems. 
\begin{theorem}[Lower bound]
\label{thm:lb}
	Let $v\in C([0,1];H^3(\R^2))\cap C^1([0,1]; L^2(\R^2))$ be a solution to~\eqref{eq:eq_v_gen} with coefficients $a_0, a_1,b_0,b_1,c_0\in L^\infty(\R^3).$ Assume that there exists a positive constant $A$ such that
	\begin{equation}\label{eq:A}
		\int_0^1 \int_{\R^2} (|v|^2 + |\partial_x v|^2 + |\partial_yv|^2 + |\Delta v|^2) \, dx dy dt\leq A^2.
	\end{equation}
	Let $\delta>0,$ $r\in (0,1/2)$ and $Q:=\{(x,y,t):\sqrt{x^2+y^2}\leq 1, t\in [r,1-r]\}$ and suppose that $\|v\|_{L^2(Q)}>\delta.$
	
	Then there exist positive constants $\widetilde{R_0},c_0,c_1$ depending on $A,$ $\|a_1\|_{L^\infty(\R^3)}, \|b_1\|_{L^\infty(\R^3)}$ and $\|c_0\|_{L^\infty(\R^3)}$ such that for $R\geq \widetilde{R_0}$
	\begin{equation}
	\label{eq:lb}
		A_R(v)\geq c_0e^{-c_1R^{3/2}}.
	\end{equation}
\end{theorem}

\begin{theorem}[Upper bound]
\label{thm:ub}
Let $v\in C([0,1];H^4(\R^2))$ be a solution of~\eqref{eq:eq_v_gen} whose coefficients $a_1,b_1$ and $c_0$ satisfy  $a_1,b_1\in L^2_xL^\infty_{yt}\cap L^1_xL^\infty_{yt}$ and $c_0\in L^\infty\cap L^2_xL^\infty_{yt},$ respectively. Then there exist positive constants $c$ and $R_0$ such that if for some $a>0$
\begin{equation*}
	v(0),v(1)\in L^2(e^{a(x^2+y^2)^{3/4}}\,dxdy),
\end{equation*}
then 
\begin{equation}
\label{eq:ub}
	A_R(v)\leq ce^{-\frac{a}{16 (28)^{3/2}}R^{3/2}}
\end{equation}
holds true for $R\geq R_0.$
\end{theorem}

Before proving Theorem~\ref{thm:lb} and Theorem~\ref{thm:ub} we show how Theorem~\ref{thm:main} follows easily as a consequence of these two results.

\begin{proof}[Proof of Theorem~\ref{thm:main}]
	We argue by contradiction. Suppose that $v:=u_1-u_2\neq 0.$ Then we can assume, after a possible translation, dilation and multiplication by a constant, that $v$ satisfies the hypotheses of Theorem~\ref{thm:lb}. Moreover, it can be seen that fixing $a_1=\frac{1}{2}(u_1 + \partial_x \mathcal{L}(u_2)),$ $b_1=\frac{1}{2}\partial_y \mathcal{L}(u_2),$ and $c_0=\frac{1}{2}\partial_x u_2,$ then $a_1,b_1\in L^2_xL^\infty_{yt}\cap L^1_xL^\infty_{yt}$ and $c_0\in L^\infty\cap L^2_xL^\infty_{yt}.$ This can be proved as in~\cite{BIM2013} (see also~\cite{CFL2019}) for the terms not involving the non-local operator $\mathcal{L}.$ For the terms involving $\mathcal{L}$ one simply uses the relation~\eqref{eq:deLRiesz} and the $L^2$-weighted boundedness~\eqref{eq:Riesz-weighted-boundedness} of the Riesz transform. Hence the hypotheses of Theorem~\ref{thm:ub} are satisfied too.  Thus, combining~\eqref{eq:lb} and~\eqref{eq:ub}, one has that for sufficiently large $R$
	\begin{equation*}
		c_0e^{c_1R^{3/2}}\leq A_R(v)\leq ce^{-\frac{a}{16(28)^{3/2}}R^{3/2}}.
	\end{equation*}
	Finally, assuming $a>a_0:=16(28)^{3/2}c_1$ and taking the limit $R$ going to infinity, we get a contradiction. Therefore $v=0$ and Theorem~\ref{thm:main} is proved.
\end{proof}

\subsection{Lower bound: Proof of Theorem~\ref{thm:lb}}
The starting point in the proof of the lower bound in Theorem~\ref{thm:lb} is a Carleman estimate for the linear operator
\begin{equation*}
	P=\partial_t + \partial_x\Delta + a_1(x,y,t) \partial_x + b_1(x,y,t) \partial_y+ c_0(x,y,t),
\end{equation*}
with $a_1,b_1, c_0\in L^\infty(\R^3).$ More precisely, we want to prove the following result.

\begin{lemma}\label{lemma:full-Carleman}
	Let $\varphi\colon [0,1]\to \R$ be a smooth function and let $D:=\R^2\times [0,1].$ We define 
	\begin{equation}\label{eq:exp-phi}
		\phi(x,y,t):=\left(\frac{x}{R} + \varphi(t) \right)^2 + \frac{y^2}{R^2}.
	\end{equation}
	Then there exists constants $c>0,$ $R_0=R_0(\|a_1\|_{L^\infty}, \|b_1\|_{L^\infty}, \|c_0\|_{L^\infty})>1$ and $\overline{C}=\max\{\|\varphi'\|_{L^\infty},\|\varphi''\|_{L^\infty},1\}>0$ such that the inequality
\begin{multline}\label{eq:full-Carl}
	\frac{\alpha^{5/2}}{R^3} \|e^{\alpha \phi}g\|_{L^2(D)}
	+ \frac{\alpha^{3/2}}{R^2}\|e^{\alpha \phi}\partial_xg\|_{L^2(D)}
	+ \frac{\alpha^{3/2}}{R^2}\|e^{\alpha \phi}\partial_yg\|_{L^2(D)}\\
	\leq c \|e^{\alpha \phi}(\partial_t + \partial_x^3 + \partial_x\partial_y^2+ a_1 \partial_x + b_1\partial_y + c_0)g\|_{L^2(D)}
\end{multline}
holds if $\alpha \geq \overline{C}R^{3/2},$ $R\geq R_0$ and $g\in C([0,1]; H^3(\R^2))\cap C^1([0,1]; L^2(\R^2))$ is compactly supported in the set $\{(x,y,t)\colon |\frac{x}{R} + \varphi(t)|\geq 1\}.$	
\end{lemma} 

As customary, the Carleman estimate contained in Lemma~\ref{lemma:full-Carleman} for the full operator $P$ can be easily proved when an analogous Carleman estimate for the leading part of $P,$ namely $\partial_t + \partial_x^3 + \partial_x\partial_y^2,$ is available. More precisely, for the proof of Lemma~\ref{lemma:full-Carleman} we will need the following result.
\begin{lemma}\label{lemma:Carleman}
	Let $\varphi\colon [0,1]\to \R$ be a smooth function and let $D\colon \R^2\times [0,1].$ We assume $R>1$ and let $\phi$ be as in~\eqref{eq:exp-phi}.
	
	Then there exist constants $c>0$ and $\overline{C}=\max\{\|\varphi'\|_{L^\infty},\|\varphi''\|_{L^\infty},1\}>0$ such that the inequality
	  \begin{equation}\label{eq:Carl}
	\frac{\alpha^{5/2}}{R^3} \|e^{\alpha \phi}g\|_{L^2(D)}
	+ \frac{\alpha^{3/2}}{R^2}\|e^{\alpha \phi}\partial_xg\|_{L^2(D)}
	+ \frac{\alpha^{3/2}}{R^2}\|e^{\alpha \phi}\partial_yg\|_{L^2(D)}\\
	\leq c\|e^{\alpha \phi}(\partial_t + \partial_x^3 + \partial_x\partial_y^2)g\|_{L^2(D)}
\end{equation}
holds if $\alpha\geq \overline{C}R^{3/2}$ and $g\in C([0,1]; H^3(\R^2))\cap C^1([0,1]; L^2(\R^2))$ is compactly supported in $\{(x,y,t)\colon |\frac{x}{R}+ \varphi(t)|\geq 1\}.$
\end{lemma} 
\begin{remark}
	A similar estimate to~\eqref{eq:Carl} was already proved in~\cite[Lemma 3.1]{BIM2013}. In their version the term on the left-hand side depending on $\partial_y$ is not present (and cannot be obtained using their proof). As we shall see below, for the purpose of proving estimate~\eqref{eq:full-Carl} in Lemma~\ref{lemma:full-Carleman} for the full operator $P$ (notice that $P$ has a lower order term involving $\partial_y$) we need the Carleman estimate in the more general form~\eqref{eq:Carl}. The proof of~\eqref{eq:Carl} slightly deviates from the analogous estimate in~\cite[Eq.~(3.2)]{BIM2013} for this reason we decided to provide the explicit proof in Appendix~\ref{ap:Carleman} below.  
\end{remark}

Now, we show how to obtain estimate~\eqref{eq:full-Carl} in Lemma~\ref{lemma:full-Carleman} from estimate~\eqref{eq:Carl} in Lemma~\ref{lemma:Carleman}
\begin{proof}[Proof of Lemma~\ref{lemma:full-Carleman}]
From estimate~\eqref{eq:Carl}, adding and subtracting the lower order terms of $P,$ namely $a_1\partial_x + b_1 \partial_y + c_0,$ and using Hölder inequality one has
\begin{multline*}
	\frac{\alpha^{5/2}}{R^3}\|e^{\alpha \phi}g\|_{L^2(D)}
	+ \frac{\alpha^{3/2}}{R^2}\|e^{\alpha\phi} \partial_xg\|_{L^2(D)}
	+ \frac{\alpha^{3/2}}{R^2}\|e^{\alpha\phi} \partial_yg\|_{L^2(D)}\\
	\lesssim_c \|e^{\alpha \phi}(\partial_t + \partial_x^3 + \partial_x \partial_y^2 + a_1 \partial_x + b_1\partial_y + c_0)g\|_{L^2(D)}\hspace{+2cm}\\
	+ \|a_1\|_{L^\infty(D)} \|e^{\alpha \phi}\partial_x g\|_{L^2(D)}
	+  \|b_1\|_{L^\infty(D)} \|e^{\alpha \phi}\partial_y g\|_{L^2(D)}
	+  \|c_0\|_{L^\infty(D)} \|e^{\alpha \phi}g\|_{L^2(D)},
\end{multline*} 
where here the symbol $\lesssim_c$ is used instead of $\leq c$ and $c$ is as in~\eqref{eq:Carl}.
Under the hypothesis $\alpha\geq \overline{C}R^{3/2},$ the ratios $\alpha^{3/2}/R^2$ and $\alpha^{5/2}/R^3$ grow as a \emph{positive} fractional power of $R,$ as a consequence, as soon as $R$ is taken sufficiently large, namely $R\geq R_0$ with $R_0$ some constant depending on the $L^\infty$-norms of $a_1, b_1$ and $c_0,$ then the last three terms on the right-hand side can be absorbed by the corresponding terms on the left-hand side giving the desired estimate.  
\end{proof}
Now we are in position to prove the lower bound in Theorem~\ref{thm:lb}.
\begin{proof}[Proof of Theorem~\ref{thm:lb}]
	As in~\cite{BIM2013} we introduce the function $\theta\in C^\infty_0(\R^2)$ such that $\theta=1$ if $\sqrt{x^2+y^2}\leq R-1,$ $\theta=0$ if $\sqrt{x^2+y^2}\geq R.$ Let $\mu \in C^\infty(\R)$ be an increasing function with $\mu=0$ in $(-\infty, 2]$ and $\mu=1$ in $[3,\infty).$ We take $\varphi\in C^\infty([0,1])$ such that $\varphi=0$ in $[0,r/2]\cup [1-r/2,1],$ $\varphi=4$ in $[r,1-r],$ $\varphi$ is increasing in $[r/2,r]$ and $\varphi$ is decreasing in $[1-r,1-r/2].$ 
	We want to apply Lemma~\ref{lemma:full-Carleman} to the function
	\begin{equation}\label{eq:def-g}
		g(x,y,t):=\theta(x,y)\mu \left( \frac{x}{R} + \varphi(t)\right) v(x,y,t),
		\quad (x,y)\in \R^2, t\in [0,1],
	\end{equation} 
	where $v$ satisfies~\eqref{eq:eq_v_gen}.
	
	It is not difficult to check that $g$ satisfies
	\begin{equation*}
			(\partial_t + \partial_x^3 + \partial_x\partial_y^2 + a_1\partial_x + b_1\partial_y + c_0)g
		=\mu \theta (\partial_t + \partial_x^3 + \partial_x\partial_y^2 + a_1\partial_x + b_1\partial_y + c_0)v + F_1 + F_2,
	\end{equation*}	
	where 
	\begin{equation*}
		F_1:=\mu \big[ 
		\partial_x(\Delta \theta) v 
		+3\partial_{x}^2 \theta \partial_x v 
		+ \partial_{y}^2\theta \partial_x v	
		+ 2\partial_{x}\partial_y\theta \partial_y v 
		+ 3\partial_x\theta \partial_{x}^2 v
		+ 2\partial_y \theta \partial_{x}\partial_{y}v
		+ \partial_x\theta \partial_{y}^2v
		+ a_1 \partial_x \theta v 
		+ b_1\partial_y \theta v	
		\big]
	\end{equation*}
	and
	\begin{multline*}
		F_2:=\big[\theta\partial_t \mu   
		+ 3 \partial_x^2 \theta \partial_x \mu
		+3 \partial_x\theta\partial_{x}^2\mu
		+ \theta \partial_{x}^3 \mu
		+\partial_y^2 \theta \partial_x \mu 
		+a_{1}\theta \partial_x \mu \big]v\\
		+[6\partial_x\theta\partial_x\mu +3 \theta\partial_x^2\mu ]\partial_x v
		+2\partial_y\theta\partial_x\mu \partial_yv
		+3\theta \partial_x \mu  \partial_x^2 v + \theta\partial_x \mu  \partial_y^2 v.		
	\end{multline*}
	Applying Lemma~\ref{lemma:full-Carleman} to $g$
	 and using that $v$ satisfies~\eqref{eq:eq_v_gen} one has
	\begin{equation}\label{eq:preliminary}
	\begin{split}
		c\frac{\alpha^{5/2}}{R^3} \|e^{\alpha \phi}g\|_{L^2(D)}
	&\leq \|e^{\alpha \phi}\theta\mu (a_0\partial_x \mathcal{L} + b_0\partial_y \mathcal{L})v \|_{L^2(D)}
	+ \|e^{\alpha \phi} F_1\|_{L^2(D)} + \|e^{\alpha \phi} F_2\|_{L^2(D)}\\
	&\leq \|a_0\|_{L^\infty(D)} \|e^{\alpha \phi}\theta\mu \partial_x \mathcal{L}v\|_{L^2(D)}
	+ \|b_0\|_{L^\infty(D)} \|e^{\alpha \phi}\theta\mu \partial_y \mathcal{L}v\|_{L^2(D)}\\  
	&\phantom{\leq}+ \|e^{\alpha \phi} F_1\|_{L^2(D)} + \|e^{\alpha \phi} F_2\|_{L^2(D)},
	\end{split}
	\end{equation}
	where, with an abuse of notation, we called $c$ the constant $c^{-1}$ and where we have used that $a_0,b_0\in L^\infty(D).$
	We want to estimate first the terms involving $F_1$ and $F_2.$ As regard with $F_1$ one observes that all the terms in $F_1$ contain derivatives of $\theta.$ This implies that $F_1$ is supported in the set $Q_R\times [0,1],$ where $Q_R$ is the anulus $Q_R:=\{(x,y)\colon R-1\leq \sqrt{x^2+y^2}\leq R\}.$ In $Q_R\times [0,1]$ one has $\phi(x,y,t)\leq 25,$ indeed
	\begin{equation*}
		\begin{split}
		\phi(x,y,t)&
		\leq \frac{x^2}{R^2} + \frac{y^2}{R^2} + \varphi(t)^2 +2|\varphi(t)|\sqrt{\frac{x^2+y^2}{R^2}}\\
		&\leq 25.
		\end{split}
	\end{equation*}
	As for $F_2,$ since the derivatives of $\mu$ appear, this is supported in the set $\{(x,y,t)\colon 2\leq |\tfrac{x}{R} + \varphi(t)|\leq 3, t\in [0,1]\}.$ From this, it follows that $\phi(x,y,t)\leq 10$ (recall that the function $\theta$ is supported in $x^2+y^2\leq R^2$). Using these facts one has 
	\begin{equation}\label{eq:F1-F2}
		\|e^{\alpha \phi} F_1\|_{L^2(D)}\leq c_1 e^{25\alpha} A_R(v),
	\quad \text{and} \quad 
		\|e^{\alpha \phi} F_2\|_{L^2(D)}\leq c_2 e^{10\alpha} A,
	\end{equation}
	where $A_R(v)$ was defined in~\eqref{eq:A_R(v)} and where $A$ is as in~\eqref{eq:A}.
	
	Now we want to estimate the terms involving the non local operators $\partial_x \mathcal{L}$ and $\partial_y \mathcal{L}.$ We will show that the operators above, being essentially the square of the Riesz transform (see~\eqref{eq:deLRiesz}), are bounded operators in a suitable weighted-$L^2$ space. More specifically, we proved the following lemma. 
	\begin{lemma}	\label{lemma:weighted-bound}
		Let $\phi$ be as in~\eqref{eq:exp-phi} and let $\theta$ and $\mu$ be as in the proof of Theorem~\ref{thm:lb}. Then the following estimates
		\begin{align}
		\label{eq:de-x}
			\|e^{\alpha \phi}\theta \mu \partial_x \mathcal{L}v\|_{L^2(D)}
			&\leq \|e^{\alpha \phi}\theta \mu v\|_{L^2(D)},\\
			\label{eq:de-y}
			\|e^{\alpha \phi}\theta \mu \partial_y \mathcal{L}v\|_{L^2(D)}
			&\leq \|e^{\alpha \phi}\theta \mu v\|_{L^2(D)}
		\end{align}
		hold.
	\end{lemma}
	In order not to weight down the presentation, we provide the proof of the lemma at the end of the section.
	Now plugging estimates~\eqref{eq:F1-F2}, ~\eqref{eq:de-x} and~\eqref{eq:de-y} in~\eqref{eq:preliminary} and recalling~\eqref{eq:def-g} gives
	\begin{equation*}
		c\frac{\alpha^{5/2}}{R^3} \|e^{\alpha \phi}g\|_{L^2(D)}
	\leq \|a_0\|_{L^\infty(D)} \|e^{\alpha \phi}g\|_{L^2(D)}
	+ \|b_0\|_{L^\infty(D)} \|e^{\alpha \phi}g\|_{L^2(D)}
	+c_1 e^{25\alpha}A_R(v) +c_2e^{10\alpha}A,
	\end{equation*}
	and since $\alpha^{5/2}/R^3$ grows as a positive power of $R,$ the first two terms in the right-hand side can be absorbed on the left-hand side giving
	\begin{equation*}
		c\frac{\alpha^{5/2}}{R^3} \|e^{\alpha \phi}g\|_{L^2(D)}
	\leq +c_1 e^{25\alpha}A_R(v) +c_2e^{10\alpha}A.
	\end{equation*}
	Since $g=v$ in $Q$ we obtain 
	\begin{equation*}
		\begin{split}
			c\frac{\alpha^{5/2}}{R^3}\|e^{\alpha \phi} g\|_{L^2(D)}
			&\geq c\frac{\alpha^{5/2}}{R^3}\|e^{\alpha \phi} g\|_{L^2(Q)}
			=c\frac{\alpha^{5/2}}{R^3}\|e^{\alpha \phi} v\|_{L^2(Q)}\\
			&\geq c\frac{\alpha^{5/2}}{R^3} e^{10\alpha}\|v\|_{L^2(Q)},
		\end{split}
	\end{equation*}
	where in the last inequality we have used that in $Q$ one has $\phi\geq 10$ for sufficiently large $R.$ From the previous two inequalities, using that $\|v\|_{L^2(Q)}\geq \delta>0$ and dividing by $e^{10\alpha}$ one has
	\begin{equation*}
		c\frac{\alpha^{5/2}}{R^3} \delta \leq c_1e^{15\alpha}A_R(v) +c_2 A.
	\end{equation*}
	Taking $\alpha=\overline{C}R^{3/2},$ with $\overline{C}$ as in Lemma~\ref{lemma:full-Carleman} one gets
	\begin{equation*}
		c \overline{C}^{5/2} R^{3/4}\delta \leq c_1e^{15\overline{C}R^{3/2}}A_R(v) +c_2 A.
	\end{equation*}
	If we take $R$ large enough, the second term of the right-hand side can be absorbed by the term on the left-hand side. Thus, we conclude that there exists $R_0>0$ such that 
	\begin{equation*}
		A_R(v) \geq Ce^{-15\overline{C}R^{3/2}}
	\end{equation*}
	holds true for $R\geq R_0.$ This concludes our proof.
	\end{proof}
	
	As already anticipated, we provide now the proof of Lemma~\ref{lemma:weighted-bound}.
	\begin{proof}[Proof of Lemma~\ref{lemma:weighted-bound}]
		 We will see explicitly how to obtain~\eqref{eq:de-x}, estimate~\eqref{eq:de-y} can be proved analogously. Taking $f:=e^{-\alpha \phi},$ $\varepsilon, \delta>0$ one has
	\begin{multline}\label{eq:eps-delta}
		\|e^{\alpha \phi}\theta \mu \partial_x \mathcal{L}v\|_{L^2(D)}\\
		\leq\|e^{\alpha \phi}(\theta + \varepsilon f) (\mu + \delta) \partial_x \mathcal{L}v\|_{L^2(D)}
		+\delta \|e^{\alpha \phi}\theta \partial_x \mathcal{L}v\|_{L^2(D)}
		+ \varepsilon\|\mu \partial_x \mathcal{L}v\|_{L^2(D)}
		+ \varepsilon \delta \|\partial_x \mathcal{L}v\|_{L^2(D)}.
	\end{multline}
	We need to estimate only the first term of the right-hand side. Indeed the remaining norms are all bounded from above by the quantity $\|\partial_x \mathcal{L}v\|_{L^2(D)}$ which, in turn, is bounded by $\|v\|_{L^2(D)}$ (this follows from~\eqref{eq:deLRiesz} and~\eqref{eq:Riesz-boundedness}). Thus the remaining terms are all negligible in the limit $\varepsilon$ and $\delta$ to zero.
	
	 Coming back to $\|e^{\alpha \phi}(\theta + \varepsilon f) (\mu + \delta) \partial_x \mathcal{L}v\|_{L^2(D)},$ recalling the relations~\eqref{eq:deLRiesz} we want to show that $w:=e^{\alpha \phi}(\theta + \varepsilon f)(\mu + \delta)\in A_2(\R^2)$ to conclude from Lemma~\ref{lemma:boundedness} the bound
	\begin{equation}\label{eq:weighted-bound}
		\|e^{\alpha \phi}(\theta + \varepsilon f) (\mu + \delta) \partial_x \mathcal{L}v\|_{L^2(D)}
		\leq \|e^{\alpha \phi}(\theta + \varepsilon f) (\mu + \delta)v\|_{L^2(D)}.
	\end{equation}
	We distinguish the case when we are inside the support of $\theta$ or outside. 
	\begin{description}
	\item[Case $\tfrac{x^2}{R^2} + \tfrac{y^2}{R^2}\leq 1$]. In this case, the following chains of inequalities hold:
	\begin{equation*}
			w=e^{\alpha \phi}\theta (\mu + \delta) + \varepsilon (\mu + \delta)\leq e^{2\alpha (1+ \|\phi\|_{L^\infty}^2)} (1+ \delta) + \varepsilon (1+ \delta),
	\end{equation*}
	and 
	\begin{equation*}
		w\geq e^{\alpha \phi} \varepsilon f (\mu + \delta) \geq \varepsilon \delta.
	\end{equation*}
	\item [Case $\tfrac{x^2}{R^2} + \tfrac{y^2}{R^2}\geq 1$]. In this case, since $\theta=0,$ one has 
		$w=\varepsilon (\mu + \delta).$
	\end{description}
With these estimates at hands we want to show that $w\in A_2(\R^2).$ In other words, we want to see that 
\begin{equation*}
	\sup_B\left(\frac{1}{|B|} \int_B w \right) \left(\frac{1}{|B|}\int_B w^{-1} \right)\leq C,
\end{equation*}
where $B$ is any ball in $\R^2$ and $C$ is a constant independent of $B$ (\emph{cfr.}~\eqref{eq:Ap}). In order to do that, for any ball $B$ in $\R^2,$ we introduce the following notations
\begin{equation*}
	B_i=B \cap \{(x,y)\colon \tfrac{x^2}{R^2} + \tfrac{y^2}{R^2}\leq 1	\},
	\qquad \text{and} \qquad
	B_e=B\cap \{(x,y)\colon \tfrac{x^2}{R^2} + \tfrac{y^2}{R^2}\geq 1	\},
\end{equation*}
therefore $B=B_i\cup B_e.$ Now
\begin{equation*}
	\begin{split}
		\left(\frac{1}{|B|} \int_B w \right) \left(\frac{1}{|B|}\int_B w^{-1} \right)
		&=\frac{1}{|B|^2} \left( \int_{B_i} w + \int_{B_e} w\right) \left( \int_{B_i} w^{-1} + \int_{B_e} w^{-1}\right)\\
		&\leq \frac{1}{|B|^2}\left(|B_i|\big(e^{2\alpha(1+ \|\phi\|_{L^\infty}^2)} + \varepsilon\big)(1+\delta) + |B_e|\varepsilon (1+\delta) \right) \left(|B_i| \frac{1}{\varepsilon \delta} + |B_e|\frac{1}{\varepsilon \delta} \right)\\
		&\leq \frac{1}{|B|^2}(c|B|)\frac{|B|}{\varepsilon \delta}\\
		&\leq c_{\varepsilon, \delta},
	\end{split}	
\end{equation*}
where $c_{\varepsilon, \delta}$ does not depend on $B.$ Using the relations~\eqref{eq:deLRiesz} and the fact that $w\in A_2(\R^2),$ applying twice Lemma~\ref{lemma:boundedness} gives~\eqref{eq:weighted-bound}.
Using estimate~\eqref{eq:weighted-bound} in~\eqref{eq:eps-delta} and letting $\varepsilon$ and $\delta$ go to zero one gets~\eqref{eq:de-x}.
	\end{proof}

\subsection{Upper bound: Proof of Theorem~\ref{thm:ub}}
This section is concerned with the proof of the upper bound in Theorem~\ref{thm:ub}. The starting point is the following persistence-decay result. 
\begin{lemma}\label{lemma:persistence-full}
	Let $w\in C([0,1]; H^4(\R^2))\cap C^1([0,1]; L^2(\R^2))$ such that for all $t\in [0,1]$ the support of $w(t)$ is contained in a compact subset $K$ of $\R^2$ and let $D:=\R^2\times [0,1].$ 
	
	Assume that $c_0\in L^\infty \cap L_x^2 L_{yt}^\infty$ and $a_1, b_1\in L_x^2L_{yt}^\infty \cap L_x^1L_{yt}^\infty,$ with small norms in these spaces.  
	Then there exists $c>0$ independent of the set $K$ such that for $\beta \geq 1$ and $\lambda \geq 7\beta$ the following estimate 
	\begin{multline}\label{eq:persistence-full}
		\|e^{\lambda |x| + \beta |y|} w\|_{L^2(D)}
		+ \sum_{0<k+l\leq 2} \|e^{\lambda |x| + \beta |y|}\partial_x^k\partial_y^l w\|_{L_x^\infty L_{yt}^2(D)}\\
		\leq c \lambda^2 \Big(\|J^3(e^{\lambda |x| + \beta |y|} w(0))\|_{L^2(\R^2)}+\|J^3(e^{\lambda |x| + \beta |y|} w(1))\|_{L^2(\R^2)} \Big)\\
		+c \|e^{\lambda |x| + \beta |y|}(\partial_t + \partial_x^3 + \partial_y \partial_y^2 + a_1\partial_x + b_1\partial_y + c_0)w\|_{L_t^1L_{xy}^2(D)\cap L_x^1L_{ty}^2(D)}	
	\end{multline}
	holds true. Here $J$ is the operator defined in Fourier space as $\mathcal{F}(Jg)(\xi,\eta):=(1 + \xi^2 + \eta^2)^{1/2}\mathcal{F}g(\xi,\eta),$ with $(\xi, \eta)$ the corresponding Fourier variables to the  variables $(x,y)$ in the physical space.
\end{lemma}
The previous result follows as a consequence of preliminary estimates involving only the leading part of the operator $P,$ namely $\partial_t + \partial_x^3 + \partial_x\partial_y^2.$ These estimates are contained in the following lemma which was proved in~\cite[Thm.1.2]{BIM2013} and which comes as a consequence of boundedness properties for the inverse of the operator $\partial_t + \partial_x^3 + \partial_y\partial_x^2.$
\begin{lemma}\label{lemma:persistence}
	Let $w$ be as in Lemma~\ref{lemma:persistence-full} and let $D:=\R^2 \times [0,1].$
	\begin{itemize}
		\item For $\lambda>0$ and $\beta>0,$
		\begin{multline}\label{eq:0-order}
			\|e^{\lambda |x| + \beta |y|} w\|_{L_t^\infty L_{xy}^2(D)}
			\leq \|e^{\lambda |x| + \beta |y|}w(0)\|_{L^2(\R^2)}
			+ \|e^{\lambda |x| + \beta |y|}w(1)\|_{L^2(\R^2)}\\
			+ \|e^{\lambda |x| + \beta |y|}(\partial_t + \partial_x^3 + \partial_y\partial_x^2)w\|_{L_t^1L_{xy}^2(D)}.
		\end{multline}
		\item There exists $c>0,$ independent of $K,$ such that for $\beta\geq 1$ and $\lambda \geq 7\beta,$
		\begin{multline}\label{eq:higher-order}
		\|e^{\lambda |x| + \beta |y|} Lw\|_{L_x^\infty L_{yt}^2(D)}
			\leq c\lambda^2 \Big(\|J^3(e^{\lambda |x| + \beta |y|}w(0))\|_{L^2(\R^2)}
			+ \|J^3(e^{\lambda |x| + \beta |y|}w(1))\|_{L^2(\R^2)}\Big)\\
			+ c\|e^{\lambda |x| + \beta |y|}(\partial_t + \partial_x^3 + \partial_y\partial_x^2)w\|_{L_x^1L_{yt}^2(D)},
		\end{multline}
		where $L$ denotes any operator in the set $\{\partial_x, \partial_y, \partial_x^2, \partial_{x y}, \partial_y^2\}.$
	\end{itemize}
\end{lemma}

Now we briefly show how to get Lemma~\ref{lemma:persistence-full} from Lemma~\ref{lemma:persistence}. 

\begin{proof}[Proof of Lemma~\ref{lemma:persistence-full}]
From estimate~\eqref{eq:0-order}, using the immediate estimates $\|\cdot\|_{L^2(D)}\leq \|\cdot\|_{L_t^\infty L_{xy}^2(D)}$ and $\|\cdot\|_{L_t^1 L_{xy}^2(D)}\leq \|\cdot\|_{L^2(D)},$ adding and subtracting the lower order part of the operator $P$ on the right-hand side and making used of H\"older's inequality, one has
\begin{multline*}
	\|e^{\lambda |x| + \beta |y|} w\|_{L^2(D)}
			\leq \|e^{\lambda |x| + \beta |y|}w(0)\|_{L^2(\R^2)}
			+ \|e^{\lambda |x| + \beta |y|}w(1)\|_{L^2(\R^2)}\\
			+ \|e^{\lambda |x| + \beta |y|}(\partial_t + \partial_x^3 + \partial_y\partial_x^2 + a_1 \partial_x + b_1 \partial_y + c_0)w\|_{L_t^1L_{xy}^2(D)}\\
			\hspace{+5cm}+\|a_1\|_{L_x^2L_{ty}^\infty(D)} \|e^{\lambda |x| + \beta |y|}\partial_x w\|_{L_x^\infty L_{yt}^2(D)}
			+\|b_1\|_{L_x^2L_{yt}^\infty(D)} \|e^{\lambda |x| + \beta |y|}\partial_y w\|_{L_x^\infty L_{yt}^2(D)}\\
			+ \|c_0\|_{L^\infty(D)} \|e^{\lambda |x| + \beta |y|}\partial_w\|_{L^2(D)}.
\end{multline*} 
Similarly, from~\eqref{eq:higher-order} one gets
\begin{multline*}
	\|e^{\lambda |x| + \beta |y|} Lw\|_{L_x^\infty L_{yt}^2(D)}
			\leq c\lambda^2 \Big(\|J^3(e^{\lambda |x| + \beta |y|}w(0))\|_{L^2(\R^2)}
			+ \|J^3(e^{\lambda |x| + \beta |y|}w(1))\|_{L^2(\R^2)}\Big)\\
			\hspace{+0.9cm}+ c\|e^{\lambda |x| + \beta |y|}(\partial_t + \partial_x^3 + \partial_y\partial_x^2 +a_1\partial_x + b_1\partial_y + c_0)w\|_{L_x^1L_{yt}^2(D)}\\
			\hspace{4.8cm}+c\|a_1\|_{L_x^1L_{yt}^\infty(D)}\|e^{\lambda |x| + \beta |y|}\partial_x w\|_{L_x^\infty L_{yt}^2(D)}
			+c\|b_1\|_{L_x^1L_{yt}^\infty(D)}\|e^{\lambda |x| + \beta |y|}\partial_y w\|_{L_x^\infty L_{yt}^2(D)}\\
			+c\|c_0\|_{L_x^2L_{yt}^\infty(D)} \|e^{\lambda |x| + \beta |y|}w\|_{L^2(D)}.
\end{multline*}
Summing the two previous estimates together we obtain
\begin{multline*}
	\|e^{\lambda |x| + \beta |y|}w\|_{L^2(D)} 
	+ \sum_{0<k+l\leq 2} \|e^{\lambda |x| + \beta |y|} \partial_x^k \partial_y^l w\|_{L_x^\infty L_{yt}^2(D)}\\
	\leq c\lambda^2 \Big(\|J^3(e^{\lambda |x| + \beta |y|}w(0))\|_{L^2(\R^2)}
			+ \|J^3(e^{\lambda |x| + \beta |y|}w(1))\|_{L^2(\R^2)}\Big)\\
			\hspace{2cm}+c\|e^{\lambda |x| + \beta |y|}(\partial_t + \partial_x^3 + \partial_y\partial_x^2 +a_1\partial_x + b_1\partial_y + c_0)w\|_{L_1^tL_{xy}^2(D)\cap L_x^1L_{yt}^2(D)}\\
			\qquad +c\|a_1\|_{L_x^2L_{yt^\infty(D)}\cap L_x^1L_{yt}^\infty(D)}\|e^{\lambda |x| + \beta |y|}\partial_x w\|_{L_x^\infty L_{yt}^2(D)}\\
			\hspace{3cm}+c\|b_1\|_{L_x^2L_{yt^\infty(D)}\cap L_x^1L_{yt}^\infty(D)}\|e^{\lambda |x| + \beta |y|}\partial_y w\|_{L_x^\infty L_{yt}^2(D)}\\
			+c\|c_0\|_{L^\infty(D)\cap L_x^2L_{yt}^\infty(D)} \|e^{\lambda |x| + \beta |y|}w\|_{L^2(D)}.\hspace{2cm}
\end{multline*}
Finally, estimate~\eqref{eq:persistence-full} follows absorbing the last three terms on the right-hand side of the previous inequality into the corresponding terms of the left-hand side using the smallness of the norms of $a_1, b_1$ and $c_0$ in the corresponding spaces.
\end{proof}

Now we are in position to prove the upper bound in Theorem~\ref{thm:ub}.
\begin{proof}[Proof of Theorem~\ref{thm:ub}]
	For $R>3$ we take $N\in \N,$ with $N>28R.$ We define $\phi_{R,N}\colon[0,\infty)\to \R$ such that $0\leq \phi_{R,N}\leq 1,$ $\phi_{R,N}=1$ in $[R+1,N]$ and $\supp \phi_{R,N}\subset (R, N+1).$ We introduce the auxiliary function
	\begin{equation*}
		w(t)(x,y):=\phi_{R,N}(\sqrt{x^2+y^2})v(t)(x,y),
	\end{equation*}
	where $v$ satisfies~\eqref{eq:eq_v_gen}.
	It is easy to see that $w$ satisfies
	\begin{equation*}
		\begin{split}
		(\partial_t + \partial_x^3 +\partial_x\partial_y^2 + a_1 \partial_x + b_1 \partial_y +c_0)w
		&=\phi_{R,N}(\partial_t + \partial_x^3 \partial_x\partial_y^2 + a_1 \partial_x + b_1 \partial_y +c_0)v + F\\
		&=-\phi_{R,N}\big( a_0\partial_x \mathcal{L}(v) +b_0\partial_y \mathcal{L}(v)\big) + F,
		\end{split}
	\end{equation*}
	where 
	\begin{multline}\label{eq:F}
		F:=\big(\partial_x^3\phi_{R,N} + \partial_x\partial_y^2\phi_{R,N} + a_1\partial_x \phi + b_1\partial_y \phi \big)v
		+\big(3\partial_x^2\phi_{R,N} + \partial_y^2\phi_{R,N}\big)\partial_x v 
		+ 2\partial_{x}\partial_y \phi_{R,N} \partial_y v\\
		+3\partial_x \phi_{R,N} \partial_x^2 v +2\partial_y \phi_{R,N} \partial_x\partial_y v
		+\partial_x\phi_{R,N} \partial_y^2v.
	\end{multline}
	Now we want to apply Lemma~\ref{lemma:persistence-full} to $w.$ In order to do that we need first to have $a_1,b_1$ and $c_0$ with small norms in the corresponding spaces (\emph{cfr.} Lemma~\ref{lemma:persistence-full}) which, in general, is not necessarily true under the generous hypotheses of Theorem~\ref{thm:ub}. For this reason we introduce an auxiliary function $\widetilde{\phi}_{R,N}$ such that $\widetilde{\phi}_{R,N}\phi_{R,N}=\phi_{R,N}$ and $\widetilde{a_1}:=a_1\widetilde{\phi}_{R,N},$ $\widetilde{b_1}:=b_1\widetilde{\phi}_{R,N}$ and $\widetilde{c_0}:=c_0\widetilde{\phi}_{R,N}$ (and also $\widetilde{a_0}:=a_0\widetilde{\phi}_{R,N}$ and $\widetilde{b_0}:=b_0\widetilde{\phi}_{R,N}$)  have small norms for $R$ sufficiently large.
	
	It is easy to show that 
	\begin{equation*}
		(\partial_t + \partial_x^3 + \partial_x\partial_y^2 + \widetilde{a_1} \partial_x + \widetilde{b_1} \partial_y +\widetilde{c_0})w
		=-\phi_{R,N}\big( \widetilde{a_0}\partial_x \mathcal{L}(v) +\widetilde{b_0}\partial_y \mathcal{L}(v)\big) + \widetilde{F},
	\end{equation*}
	where $\widetilde{F}$ has the same form of $F$ but with $a_1$ and $b_1$ replaced by $\widetilde{a_1}$ and $\widetilde{b_1}.$
	
	Now we can apply Lemma~\ref{lemma:persistence-full} obtaining that for $\beta\geq 1$ and $\lambda\geq 7\beta$ it holds
	\begin{multline}\label{eq:pers-preliminary}
	\|e^{\lambda |x| + \beta |y|} w\|_{L^2(D)}
		+ \sum_{0<k+l\leq 2} \|e^{\lambda |x| + \beta |y|}\partial_x^k\partial_y^l w\|_{L_x^\infty L_{yt}^2(D)}\\
		\leq c \lambda^2 \Big(\|J^3(e^{\lambda |x| + \beta |y|} w(0))\|_{L^2(\R^2)}+\|J^3(e^{\lambda |x| + \beta |y|} w(1))\|_{L^2(\R^2)} \Big)\\
		\hspace{1cm}+c\|\widetilde{a_0}\|_{L^\infty(D)\cap L_x^2L_{yt}^\infty(D)} \|e^{\lambda |x| + \beta |y|}\phi_{R,N}\partial_x\mathcal{L}(v)\|_{L^2(D)}\\
		\hspace{1cm}+ c\|\widetilde{b_0}\|_{L^\infty(D)\cap L_x^2L_{yt}^\infty(D)} \|e^{\lambda |x| + \beta |y|}\phi_{R,N}\partial_y\mathcal{L}(v)\|_{L^2(D)}\\
		+c \|e^{\lambda |x| + \beta |y|}\widetilde{F}\|_{L_t^1L_{xy}^2(D)\cap L_x^1L_{ty}^2(D)},	\hspace{4cm}
	\end{multline}
	where we have used that by H\"older's inequality 
	\begin{multline*}
		\|e^{\lambda |x| + \beta |y|}\phi_{R,N}(\widetilde{a_0}\partial_x\mathcal{L}(v) +\widetilde{b_0}\partial_y\mathcal{L}(v))\|_{L_t^1L_{xy}^2(D)\cap L_x^1L_{ty}^2(D)}\\
		\leq\|\widetilde{a_0}\|_{L^\infty(D)\cap L_x^2L_{yt}^\infty(D)} \|e^{\lambda |x| + \beta |y|}\phi_{R,N}\partial_x\mathcal{L}(v)\|_{L^2(D)}
		+ \|\widetilde{b_0}\|_{L^\infty(D)\cap L_x^2L_{yt}^\infty(D)} \|e^{\lambda |x| + \beta |y|}\phi_{R,N}\partial_y\mathcal{L}(v)\|_{L^2(D)}.
	\end{multline*}
	We start considering the terms involving the non-local operators $\partial_x \mathcal{L}$ and $\partial_y\mathcal{L}.$ We will see explicitly how to treat $\|e^{\lambda |x| + \beta |y|} \phi_{R,N} \partial_x\mathcal{L}(v)\|_{L^2(D)},$ the corresponding term in the $y$ variable can be estimated analogously. As in the proof of Lemma~\ref{lemma:weighted-bound} above, we take an auxiliary function  $f:=e^{-(\lambda |x| + \beta |y|)}$ and $\varepsilon>0$ and we write
	\begin{equation}\label{eq:eps-perturbation}
		\|e^{\lambda |x| + \beta |y|} \phi_{R,N} \partial_x\mathcal{L}(v)\|_{L^2(D)}
		\leq \|e^{\lambda |x| + \beta |y|} (\phi_{R,N} + \varepsilon f) \partial_x\mathcal{L}(v)\|_{L^2(D)}  
		+ \varepsilon \|\partial_x\mathcal{L}(v)\|_{L^2(D)}.
	\end{equation}
	Using $L^2$-boundedness of the Riesz transform (see~\eqref{eq:Riesz-boundedness}) one has that the last term on the right-hand side of~\eqref{eq:eps-perturbation} is finite and thus negligible due to the arbitrariness of $\varepsilon>0.$ For estimating the first term on the right-hand side of~\eqref{eq:eps-perturbation} we want to show that the weight $w:=e^{\lambda |x| + \beta |y|} (\phi_{R,N} + \varepsilon f)\in A_2(\R^2)$ to conclude from the bound~\eqref{eq:Riesz-weighted-boundedness} in Lemma~\ref{lemma:boundedness} that 
	\begin{equation}\label{eq:non-local-pers}
		\|e^{\lambda |x| + \beta |y|} (\phi_{R,N} + \varepsilon f) \partial_x\mathcal{L}(v)\|_{L^2(D)}\leq \|e^{\lambda |x| + \beta |y|} (\phi_{R,N} + \varepsilon f) v\|_{L^2(D)}.  
	\end{equation}
	To show that, as in the proof of Lemma~\ref{lemma:weighted-bound}, we distinguish the case when we are inside the support of $\phi_{R,N}$ or outside.
	\begin{description}
	\item[Case] $\sqrt{x^2+y^2}\in [R,N+1].$ In this case, the following chains of inequalities hold:
	\begin{equation*}
		w=e^{\lambda |x| + \beta |y|} (\phi_{R,N} + \varepsilon f)
		\leq e^{(\lambda + \beta)(N+1)} +\varepsilon
	\end{equation*}
	and 
	\begin{equation*}
		w\geq \varepsilon.
	\end{equation*}
	\item[Case] $\sqrt{x^2+y^2}\notin [R, N+1].$ In this case, since $\phi_{R,N}=0,$ one has $w=\varepsilon.$
	\end{description}
	These bounds show that $w\in A_2(\R^2)$ (see Lemma~\ref{lemma:weighted-bound} above for the details), thus estimate~\eqref{eq:non-local-pers} holds. Since $\varepsilon$ is arbitrary, letting $\varepsilon$ go to zero, one has from~\eqref{eq:non-local-pers} and~\eqref{eq:eps-perturbation}
	\begin{equation*}
		\|e^{\lambda |x| + \beta |y|} \phi_{R,N} \partial_x\mathcal{L}(v)\|_{L^2(D)}\leq \|e^{\lambda |x| + \beta |y|} w\|_{L^2(D)}.
	\end{equation*} 
	Using this bound and the fact that $\widetilde{a_0}$ and $\widetilde{b_0}$ have small norms, the non-local terms on the right-hand side of~\eqref{eq:pers-preliminary} can be absorbed on the left-hand side giving
	\begin{multline}\label{eq:first-ub}
	\|e^{\lambda |x| + \beta |y|} w\|_{L^2(D)}
		+ \sum_{0<k+l\leq 2} \|e^{\lambda |x| + \beta |y|}\partial_x^k\partial_y^l w\|_{L_x^\infty L_{yt}^2(D)}\\
		\leq c \lambda^2 \Big(\|J^3(e^{\lambda |x| + \beta |y|} w(0))\|_{L^2(\R^2)}+\|J^3(e^{\lambda |x| + \beta |y|} w(1))\|_{L^2(\R^2)} \Big)\\
		+c \|e^{\lambda |x| + \beta |y|}\widetilde{F}\|_{L_t^1L_{xy}^2(D)\cap L_x^1L_{ty}^2(D)}	\hspace{4cm},
	\end{multline}

We continue by estimating the term depending on $\widetilde{F}.$ From the explicit expression of $\widetilde{F}$ (see~\eqref{eq:F}) one notice that only derivatives of $\phi_{R,N}$ appear. Thus the support of $\widetilde{F}$ is the union of the sets $\Omega_R:=\{(x,y)\colon \sqrt{x^2+y^2}\in (R, R+1)\}$ and $\Omega_N:=\{(x,y)\colon\sqrt{x^2+y^2}\in (N,N+1)\},$ whose areas are of order $R$ and $N$ respectively. In the first set $\Omega_R$ one has
\begin{equation*}
	\lambda|x|+ \beta|y|
	\leq(\lambda + \beta)\sqrt{x^2+y^2}
	\leq (\lambda +\beta)(R+1).
\end{equation*} 
From this estimate, using Cauchy-Schwarz inequality,  one gets
\begin{equation}\label{eq:omegaR}
	\begin{split}
	\|e^{\lambda |x| + \beta |y|}\widetilde{F} \chi_{\Omega_R}\|_{L_t^1L_{xy}^2(D)\cap L_x^1L_{ty}^2(D)}
	&\leq e^{(\lambda +\beta)(R+1)} \|\widetilde{F} \chi_{\Omega_R}\|_{L_t^1L_{xy}^2(D)\cap L_x^1L_{ty}^2(D)}\\
	&\leq c R^{1/2} e^{(\lambda +\beta)(R+1)}\||v| + |\nabla v| + |\Delta v| + |\partial_x\partial_y v|\|_{L^2(D)}\\
	&\leq cR^{1/2} e^{(\lambda +\beta)(R+1)}, 
	\end{split}
\end{equation}
here the constant $c>0$ changes from line to line.
For the estimate in $\Omega_N$ one observe that since $v(0), v(1)\in L^2(e^{a(x^2+y^2)^{3/4}}dxdy)$ it follows that for all $\overline{\lambda}>0$ and all $\overline{\beta}>0,$ $v(0), v(1)\in L^2(e^{2\overline{\lambda}|x| + 2\overline{\beta} |y|}dx dy).$ Using a similar argument to the one for the classical ZK contained in~\cite[Thm. 1.3]{BIM2011}, one can prove that this property is preserved in the whole interval of time $[0,1]$. More precisely, one has that $v$ is a bounded function from $[0,1]$ with values in $H^3(e^{2\overline{\lambda}|x| + 2\overline{\beta}|y|}dx dy).$ Thus, we can take $\overline{\lambda}=\lambda+ 1$ and $\overline{\beta}=\beta +1$ and applying the Cauchy-Schwarz inequality one obtains
\begin{equation}\label{eq:omegaN}
	\|e^{\lambda |x| + \beta |y|}\widetilde{F} \chi_{\Omega_N}\|_{L_t^1L_{xy}^2(D)\cap L_x^1L_{ty}^2(D)}
	\leq C_{\lambda,\beta} N^{1/2} e^{-N}\|e^{(\lambda+1) |x| + (\beta+1) |y|}\widetilde{F}\|_{L^2(D)},
\end{equation}
where the last term is finite. This implies that, since $N$ can be arbitrary large, the right-hand side of~\eqref{eq:omegaN} is negligible. Thus, taking the limit as $N$ tends to infinity, from~\eqref{eq:omegaR} and~\eqref{eq:omegaN} one gets
\begin{equation}\label{eq:est-F}
 	\|e^{\lambda |x| + \beta |y|}\widetilde{F}\|_{L_t^1L_{xy}^2(D)\cap L_x^1L_{ty}^2(D)} \leq C R^{1/2}e^{(\lambda +\beta)(R+1)}.
\end{equation}
Before providing the bound for the terms depending on the initial, respectively, final data $w(0)$ and $w(1)$ on the right-hand side of~\eqref{eq:first-ub} we make the following observation: considering the region $Q_{28R}:=\{(x,y,t)\colon 28R-1\leq \sqrt{x^2+y^2}\leq 28R\}.$ One has $Q_{28R}\times [0,1]\subset D_R,$ moreover, since $N>28R,$ $v=w$ in $Q_{28R}\times [0,1].$ In $Q_{28R}\times [0,1]$ we also have  the following bound
\begin{equation*}
	\lambda|x| + \beta|y|\geq \frac{\lambda}{7}(|x|+|y|)
	\geq \frac{\lambda}{7}(28R-1)\geq 3\lambda R.
\end{equation*}		 
From that, employing the H\"older inequality, the bound~\eqref{eq:first-ub} and the estimate~\eqref{eq:est-F}, one has
\begin{equation}\label{eq:last-one}
\begin{split}
	e^{3\lambda R}A_{28R}
	&\leq \|e^{\lambda|x|+ \beta |y|} v\|_{L^2(Q_{28}\times [0,1])}
	+ \sum_{0<k+l\leq 2} \|e^{\lambda|x|+ \beta |y|} \partial_x^k \partial_y^lv\|_{L^2(Q_{28}\times [0,1])}\\
	&\leq R^{1/2} \Big( \|e^{\lambda|x|+ \beta |y|} w\|_{L^2(D)}
	+ \sum_{0<k+l\leq 2} \|e^{\lambda|x|+ \beta |y|} \partial_x^k \partial_y^lw\|_{L^2(D)}\Big)\\
	&\leq c \lambda^2 R^{1/2} \Big(\|J^3(e^{\lambda |x| + \beta |y|} w(0))\|_{L^2(\R^2)}+\|J^3(e^{\lambda |x| + \beta |y|} w(1))\|_{L^2(\R^2)} \Big)
		+ce^{2\lambda R},
\end{split}
\end{equation}
where in the last inequality we have used that choosing $\beta=\lambda/7,$ $\lambda \geq 7$ and using that $R>3$ one has
\begin{equation*}
	Re^{(\lambda + \beta) (R+1)}
	\leq ce^{(\lambda+\beta +1)(R+1)}
	\leq ce^{\lambda\big(1+\tfrac{1}{7} +\tfrac{1}{7}\big)\tfrac{4}{3}R}
	\leq ce^{2\lambda R}.
\end{equation*}

	 It remains to estimate the terms depending on the initial, respectively, final data $w(0)$ and $w(1).$ We shall see explicitly only the one for $w(0)$ as the one for $w(1)$ can be treated analogously.
	We will make use of the following interpolation result that we proved at the end of this section (see also~\cite[Lemma 1]{FP2011} for a similar result).
\begin{lemma}\label{lemma:interpolation}
	For $k\in \N$ and $\beta>0$ let $f\in H^k(\R^2)\cap L^2(e^{2\beta (|x|+ |y|)}dx dy).$ Then for $\theta\in [0,1]$ one has
	\begin{equation}\label{eq:interpolation}
		\|J^{\theta k}(e^{(1-\theta)\beta (|x|+|y|)}f)\|_{L^2(\R^2)}
		\leq c\|J^k f\|_{L^2(\R^2)}^\theta \|e^{\beta(|x|+|y|)} f\|_{L^2(\R^2)}^{1-\theta},
	\end{equation}
	with a constant $c>0$ depending on $k$ and $\beta.$
	\end{lemma}		
	Using the previous result for $k=4,$ $\beta=4\lambda$ and $f=w(0)$ one has
	\begin{equation*}
		\|J^3(e^{\lambda (|x|+|y|)}w(0))\|_{L^2(\R^2)}
		\leq c\|J^4 w(0)\|_{L^2(\R^2)}^{3/4} \|e^{4\lambda(|x|+|y|)} w(0)\|_{L^2(\R^2)}^{1/4}.
	\end{equation*}
	First of all,  observe that since $\phi_{R,N}$ and its derivatives are bounded by a constant independent of $R$ and $N$ one has
	\begin{equation}\label{eq:inter-w0}
		\|J^4 w(0)\|_{L^2(\R^2)}\leq c \|J^4 v(0)\|_{L^2(\R^2)},
	\end{equation}
	thus  $w(0)\in H^4(\R^2)$ being $v(0)\in H^4(\R^2)$ by hypotheses.
	
	Now, since $w$ is supported  in the set
	\begin{equation*}
		D_R:=\{(x,y)\colon x^2+y^2\geq R^2\}\times [0,1],
	\end{equation*} 
	and using again that $\phi_{R,N}$ is bounded by a constant independent of $R$ and $N$ one has
	\begin{equation}\label{eq:H4}
		\|e^{4\lambda(|x|+|y|)} w(0)\|_{L^2(\R^2)}
		=\|e^{4\lambda(|x|+|y|)} w(0)\|_{L^2(x^2+y^2\geq R^2)}
		\leq c \|e^{4\lambda(|x|+|y|)} v(0)\|_{L^2(x^2+y^2\geq R^2)}.
	\end{equation}
	Moreover, choosing $\lambda:=\frac{a}{16}R^{1/2},$ if $x^2+y^2\geq R^2$ one gets
	\begin{equation}\label{eq:leq}
		4\lambda(|x| +|y|)\leq \frac{1}{2\sqrt{2}} aR^{1/2}\sqrt{x^2+y^2}
		\leq \frac{1}{2\sqrt{2}} a (x^2+y^2)^{3/4}.
	\end{equation}
	From~\eqref{eq:inter-w0}, using~\eqref{eq:H4} and~\eqref{eq:leq} one gets
	\begin{equation}\label{eq:est-w0-w1}
	\begin{split}
		\lambda^2R^{1/2}\|J^3(e^{\lambda (|x|+|y|)}w(0))\|_{L^2(\R^2)}
		&\leq c\lambda^2R^{1/2}
		\|e^{4\lambda(|x|+|y|)} w(0)\|_{L^2(\R^2)}^{1/4}\\
		&\leq c a^2 R^3\Big \|e^{\tfrac{a}{2\sqrt{2}} (x^2+y^2)^{3/4}}v(0)\Big \|_{L^2(x^2+y^2\geq R^2)}^{1/4}\\
		&\leq c a^2\big\|e^{\tfrac{a}{2} (x^2+y^2)^{3/4}}v(0)\big\|_{L^2(D)}^{1/4}\\
		&\leq c_a,
	\end{split}	
	\end{equation}
	where in the last inequality we have used that $v(0)\in L^{2}(e^{a(x^2+y^2)^{3/4}}dxdy).$ 
	Using~\eqref{eq:est-w0-w1} in~\eqref{eq:last-one} gives
	\begin{equation*}
		e^{3\lambda R} A_{28R}\leq c_a + ce^{2\lambda R}\leq c_a e^{2\lambda R}.
	\end{equation*}
	This gives the upper bound
	\begin{equation*}
		A_{28R}\leq c_a e^{-\frac{a}{16}R^{3/2}},
	\end{equation*}
	which was what we wanted to prove.
\end{proof}

As anticipated, we want now to prove the interpolation Lemma~\ref{lemma:interpolation}.
\begin{proof}[Proof of Lemma~\ref{lemma:interpolation}]
It will suffice to consider the case $k=1.$ For $\theta\in [0,1]$ and $\beta>0$ we define $w_\theta(x,y):=e^{(1-\theta)\beta(|x|+|y|)}.$ To prove inequality~\eqref{eq:interpolation} we only need to show that 
\begin{equation}\label{eq:comm-estimate}
	\|J^{\theta}(w_\theta f)\|_{L^2(\R^2)}\simeq \|w_\theta J^{\theta}f\|_{L^2(\R^2)}. 
\end{equation}
Indeed, as soon as~\eqref{eq:comm-estimate} is proved, inequality~\eqref{eq:interpolation} follows using standard interpolation of weighted Sobolev spaces. Observe that to prove~\eqref{eq:comm-estimate} it is enough to show the inequality
\begin{equation}\label{eq:comm-ineq}
	\|J^{\theta}(w_\theta f)\|_{L^2(\R^2)}\lesssim \|w_\theta J^{\theta}f\|_{L^2(\R^2)}
\end{equation}
as the reverse inequality is obtained by duality.
The inequality~\eqref{eq:comm-ineq} is equivalent to 
\begin{equation}\label{eq:equivalent-comm}
\|J^{\theta}w_\theta J^{-\theta} w_\theta^{-1}f\|_{L^2(\R^2)}\lesssim \|f\|_{L^2(\R^2)}.
\end{equation}
To prove~\eqref{eq:equivalent-comm}, we define the analytic family of operators
$T_z:=J^{z}w_z J^{-z} w_z^{-1},$ for $\Re z\in [0,1].$
Hence estimate~\eqref{eq:equivalent-comm} follows from the Stein-Weiss interpolation theorem as soon as one proves $L^2-L^2$ boundedness of $T_{i\eta}$ and $T_{1+i\eta},$ $\eta \in \R.$ In order to do that we will consistently exploit the fact that, for any $\gamma \in \R,$ the operator $J^{i\gamma}$ is bounded in $L^2$ and  it is bounded in the $w$-weighted $L^2$ space for any weight $w\in L^1_{\text{loc}},$ namely it holds true the following estimate
\begin{equation*}
	\|w J^{i\gamma}f\|_{L^2}\leq c(\gamma) \|wf\|_{L^2},
\end{equation*}
or equivalently
\begin{equation*}
	\|w f\|_{L^2}\leq c \|wJ^{-i\gamma}f\|_{L^2},
\end{equation*}
with $c$ independent of $\gamma.$

We start with the estimate for 
	$T_{i\eta}=J^{i\eta }e^{(1-i\eta)\beta(|x|+|y|)}J^{-i\eta} e^{(-1+i\eta)\beta(|x|+|y|)}.$
Using the boundedness properties of $J^{i\gamma}$ one has
\begin{equation*}
	\begin{split}
		\|J^{i\eta }e^{(1-i\eta)\beta (|x|+|y|)}f\|_{L^2(\R^2)}
		&\lesssim \|e^{(1-i\eta)\beta (|x|+|y|)}f\|_{L^2(\R^2)}\\
		&\lesssim \|e^{(1-i\eta)\beta (|x|+|y|)}J^{i\eta}f\|_{L^2(\R^2)},
	\end{split}
\end{equation*}
or equivalently
\begin{equation*}
\|T_{i\eta}f\|_{L^2(\R^2)}\lesssim \|f\|_{L^2(\R^2)}.
\end{equation*}
We consider now the operator
	$T_{1+i\eta}=J^{(1+i\eta) }e^{i\eta\beta(|x|+|y|)}J^{-(1+i\eta) } e^{-i\eta\beta(|x|+|y|)}.$
A similar argument used for the estimate of $T_{i\eta}$ gives
\begin{equation}\label{eq:J1plus}
\begin{split}
	\|J^{(1+i\eta)} e^{i\eta \beta (|x|+|y|)}f\|_{L^2(\R^2)}
	&\lesssim \|J e^{i\eta  \beta (x+y)}f\|_{L^2(x\geq 0, y\geq 0)}
	+ \|J e^{i\eta  \beta (x-y)}f\|_{L^2(x\geq 0, y< 0)}\\
	&\phantom{\quad}+ \|J e^{i\eta  \beta (-x+y)}f\|_{L^2(x< 0, y\geq 0)}
	+ \|J e^{i\eta  \beta (-x-y)}f\|_{L^2(x< 0, y< 0)}\\
	&\lesssim \|J e^{i\eta  \beta (x+y)}f\|_{L^2(\R^2)}
	+ \|J e^{i\eta  \beta (x-y)}f\|_{L^2(\R^2)}\\
	&\phantom{\quad}+ \|J e^{i\eta  \beta (-x+y)}f\|_{L^2(\R^2)}
	+ \|J e^{i\eta  \beta (-x-y)}f\|_{L^2(\R^2)}\\
	&:=\emph{r-h-s}.
\end{split}
\end{equation}
Recall that, by definition, for $s\in \R$ we have the following chain of identities 
\begin{equation*}
	\|J^sf\|_{L^2(\R^2)}=\|(1-\Delta)^{s/2}f\|_{L^2(\R^2)}\simeq \|f\|_{L^2(\R^2)} + \|D^sf\|_{L^2(\R^2)},
	\qquad D^s:=(-\Delta)^{s/2}.
\end{equation*}
Using this fact for $s=1$ and noticing that 
	$|D e^{i\eta \beta(\pm x\pm y)}|\lesssim |\eta| ,$
the Leibniz rule gives
\begin{equation}\label{eq:rhs}
	\begin{split}
		\emph{r-h-s}&\lesssim |\eta|\|f\|_{L^2(\R^2)}+ \|Df\|_{L^2(\R^2)}\\
		&\lesssim (1+ |\eta|)\|Jf\|_{L^2}\\
		&\lesssim (1+ |\eta|)\|e^{i\eta\beta (|x|+ |y|)}J^{(1+i\eta)}f\|_{L^2}.
	\end{split}
\end{equation}
Plugging~\eqref{eq:rhs} in~\eqref{eq:J1plus} one has
\begin{equation*}
	\|T_{1+i\eta}f\|_{L^2(\R^2)}\lesssim (1+|\eta|) \|f\|_{L^2(\R^2)},
\end{equation*}
which concludes the proof.
\end{proof}

\appendix
\section*{Carleman estimate: Proof of Lemma~\ref{lemma:Carleman}}
\label{ap:Carleman}
To conclude, we give the proof of the Carleman estimate contained in Lemma~\ref{lemma:Carleman}. As already mentioned, the proof is an adaptation of Lemma 3.1 in~\cite{BIM2013}. 
\begin{proof}[Proof of Lemma~\ref{lemma:Carleman}]
For notation simplicity, for the moment we make it no longer explicit the dependence on $\alpha$ of the exponential weight $e^{\alpha \phi}.$ To do that we define $\psi(x,y,t):=\alpha \phi(x,y,t).$
 
As customary, to prove~\eqref{eq:Carl} it is convenient to consider the conjugated operator
\begin{equation}\label{eq:conj-oper}
	e^{\psi(x,y,t)} (\partial_t + \partial_x^3 + \partial_x \partial_y^2) e^{-\psi(x,y,t)}.
\end{equation}
We define $f:=e^{\psi}g.$ Since $e^{\psi} \partial_x e^{-\psi}=(\partial_x - \psi_x)f$ and $e^{\psi} \partial_y e^{-\psi}=(\partial_y - \psi_y)f,$ to establish~\eqref{eq:Carl} it is sufficient to prove the corresponding estimate for the conjugated operator~\eqref{eq:conj-oper}, namely
\begin{equation}\label{eq:Carl-conj}
	\frac{\alpha^{5/2}}{R^3}\|f\|_{L^2(D)}
	+ \frac{\alpha^{3/2}}{R^2}\|(\partial_x - \psi_x)f\|_{L^2(D)}
	+ \frac{\alpha^{3/2}}{R^2}\|(\partial_y - \psi_y)f\|_{L^2(D)}\\
	\leq c \|e^{\psi}(\partial_t + \partial_x^3 + \partial_x \partial_y^2)e^{-\psi}f\|_{L^2(D)}.
\end{equation}

First of all, a straightforward computation gives
\begin{equation*}
	\begin{split}
		e^\psi (\partial_t + \partial_x^3 + \partial_x\partial_y^2)e^{-\psi}f
		&=(\partial_t -\psi_t)f 
		+(\partial_x-\psi_x)^3f
		+ (\partial_x -\psi_x)(\partial_y-\psi_y)^2f\\
		&=\partial_t f -\psi_t f
		+ \partial_x^3 f-3\psi_{x x} \partial_xf
		-3\psi_x \partial_x^2 f + 3\psi_x \psi_{xx} f 
		+ 3\psi_x^2 \partial_x f
		-\psi_x^3 f\\
		&\phantom{=}
		+\partial_x\partial_y^2 f -\psi_{yy} \partial_x f
		-2\psi_y \partial_x\partial_y f
		+\psi_y^2 \partial_x f
		-\psi_x \partial_y^2 f 
		+\psi_x \psi_{yy} f
		+2\psi_x\psi_y \partial_y f
		-\psi_x \psi_y^2 f,
	\end{split}
\end{equation*} 
where here we have used that $\psi_{xxx}=\psi_{xxy}=\psi_{xy}=0.$
We can write this as 
\begin{equation*}
	e^\psi (\partial_t + \partial_x^3 + \partial_x\partial_y^2)e^{-\psi}f 
	=\mathcal{A}f + \mathcal{S}f,
\end{equation*}
where $\mathcal{A}f$ and $\mathcal{S}f$ are respectively the anti-symmetric and symmetric operators given by
\begin{align*}
	\mathcal{A}&:=\partial_x^3 + 3 \psi_x^2 \partial_x + 3\psi_x \psi_{xx} + \partial_t +\psi_{yy}\psi_x + \psi_y^2 \partial_x + 2\psi_y \psi_x \partial_y + \partial_y^2 \partial_x\\
	\mathcal{S}&:=-3\partial_x(\psi_x \partial_x) -\psi_x^3 -\psi_t -\psi_y^2\psi_x -\psi_{yy}\partial_x -2\psi_y\partial_x\partial_y -\psi_x \partial_y^2.
\end{align*}

The main point in the proof of~\eqref{eq:Carl-conj} is to provide a lower bound for the positivity of the commutator between the symmetric and the anti-symmetric part of the conjugated operator~\eqref{eq:conj-oper}.  Indeed, it is easy to see that 
\begin{equation}\label{eq:comm}
	\begin{split}
		\|e^\psi(\partial_t + \partial_x^3 + \partial_x\partial_y^2)e^{-\psi}f\|_{L^2}^2
		&=\|(\mathcal{A} + \mathcal{S})f\|_{L^2}^2\\
		&=\|\mathcal{A}f\|_{L^2}^2 + \|\mathcal{S}f\|_{L^2}^2
		+\langle [\mathcal{S}, \mathcal{A}]f,f \rangle\\
		&\geq \langle [\mathcal{S}, \mathcal{A}]f,f \rangle, 
	\end{split}
\end{equation}	
where $\langle \cdot, \cdot\rangle$ denotes the inner product in $L^2(D).$

An explicit computation of the commutator in~\eqref{eq:comm} gives
\begin{equation}\label{eq:preliminary-comm}
	\begin{split}
		\langle [\mathcal{S}, \mathcal{A}]f,f \rangle
		=&\int_D (9\psi_x^4 \psi_{xx} -3\psi_{xx}^3 + 6\psi_{xt}\psi_{x}^2
		+\psi_{tt} + 6 \psi_x^2 \psi_y^2 \psi_{xx} + 2 \psi_{xt}\psi_y^2\\
		&\phantom{\int_D (} 
		+4\psi_x^2\psi_y^2\psi_{yy} + \psi_y^4\psi_{xx} +\psi_{xx}\psi_{yy}^2 -6 \psi_{xx}^2\psi_{yy}
		)f^2\\
		&+\int_D (18\psi_x^2\psi_{xx}-6\psi_{xt}-6\psi_y^2\psi_{xx} + 4 \psi_{yy}\psi_{y}^2) (\partial_x f)^2\\
		&+\int_D (2\psi_y^2 \psi_{xx} -6\psi_x^2 \psi_{xx}-2\psi_{xt} + 4\psi_x^2\psi_{yy})(\partial_y f)^2\\
		&+ \int_D 24\psi_x\psi_y\psi_{xx} \partial_x f \partial_y f\\
		&+\int_D (4\psi_{yy} + 6\psi_{xx})(\partial_y\partial_x f)^2
		+\int_D 9\psi_{xx}(\partial_x^2 f)^2 
		+ \int_D \psi_{xx}(\partial_y^2 f)^2.
	\end{split}
\end{equation}
We consider first the term involving $\partial_x f \partial_y f$.   Since $\psi_{xx}$ and $\psi_{yy}$ are constants,  integrating by parts yields
\begin{equation*}
	\begin{split}
		24\psi_{xx} \int_{D} \psi_x \psi_y \partial_x f \partial_y f
		&=16 \psi_{xx} \int_{D} \psi_x \psi_y \partial_x f \partial_y f
		+ 8\psi_{xx} \int_{D} \psi_x \psi_y \partial_x f \partial_y f\\
		&=8\psi_{xx}^2\psi_{yy} \int_{D}f^2
		-16 \psi_{xx} \int_{D} \psi_x \psi_y f \partial_x \partial_y f 
		+ 8\psi_{xx} \int_D \psi_x \psi_y \partial_x f \partial_y f.
	\end{split}
\end{equation*}
Observe that in the last identity, the first term is positive, instead the remaining two terms, in principle, do not have a sign. We also notice that the integrands in the remaining integrals can be thought as a double product in the development of a square. We are going to explore this idea in the following.

Since $\psi_{xx}=\psi_{yy}$ and they are positive constants, one has  
\begin{equation*}
	\begin{split}
	\int_D \big[(6 \psi_{xx} + 4\psi_{yy})\psi_x^2 \psi_y^2f^2
	-16 \psi_{xx} &\psi_x \psi_y f \partial_x \partial_y f
	+ (4\psi_{xx} + 6 \psi_{yy}) (\partial_y \partial_x f)^2\big]\\
	&=\psi_{xx} \int_{D} \big[ 10 \psi_x^2\psi_y^2 f^2 -16\psi_x\psi_y f \partial_x\partial_y f + 10 (\partial_x\partial_y f)^2\big]\\
	&=\psi_{xx} \int_{D} \left( \frac{8}{3}\psi_x\psi_y f -3\partial_x\partial_y f \right)^2 
	+\frac{26}{9} \psi_{xx} \int_D  \psi_x^2 \psi_y^2 f^2
	+ \psi_{xx} \int_D (\partial_x\partial_y f)^2\\
	&\geq \frac{26}{9} \psi_{xx} \int_D \psi_x^2 \psi_y^2 f^2
	+ \psi_{xx} \int_D (\partial_x\partial_y f)^2.
	\end{split}
\end{equation*}
Now we consider
\begin{equation*}
	\begin{split}
	\psi_{xx}\int_{D} \big[ 18\psi_x^2 (\partial_x f)^2 + 2\psi_y^2 (\partial_y f)^2 &+ 8 \psi_x \psi_y \partial_x f \partial_y f \big]\\
	&=\psi_{xx} \int_{D} (4\psi_x \partial_x f + \psi_y \partial_y f)^2
	+ 2\psi_{xx} \int_{D} \psi_x^2 (\partial_x f)^2 
	+ \psi_{xx} \int \psi_y^2 (\partial_y f)^2\\
	&\geq 2\psi_{xx} \int_{D} \psi_x^2 (\partial_x f)^2 
	+ \psi_{xx} \int \psi_y^2 (\partial_y f)^2.  
	\end{split}
\end{equation*}
Plugging the previous estimates in~\eqref{eq:preliminary-comm} one gets
\begin{equation*}
\begin{split}
		\langle [\mathcal{S}, \mathcal{A}]f,f \rangle
		 \ge &\int_D (9\psi_x^4 \psi_{xx} + 6\psi_{xt}\psi_{x}^2
		+\psi_{tt} + \tfrac{26}{9} \psi_x^2 \psi_y^2 \psi_{xx} + 2 \psi_{xt}\psi_y^2
		+ \psi_y^4\psi_{xx}
		)f^2\\
		&+\int_D (2\psi_x^2\psi_{xx}-6\psi_{xt}-2\psi_y^2\psi_{xx}) (\partial_x f)^2\\
		&+\int_D (\psi_y^2 \psi_{xx} -2\psi_x^2 \psi_{xx}-2\psi_{xt})(\partial_y f)^2\\
		&+ \psi_{xx}\int_D (\partial_y\partial_x f)^2
		+9\psi_{xx}\int_D (\partial_x^2 f)^2 
		+\psi_{xx} \int_D (\partial_y^2 f)^2.
	\end{split}
\end{equation*}
Using that $\int\psi_y^2 (\partial_x f)^2=-\int \psi_y^2 f \partial_x^2 f$ one has
\begin{equation*}
	\begin{split}
	\psi_{xx}\int_{D} \big[ \psi_y^4 f^2 -2\psi_y^2 (\partial_x f)^2 + 9 (\partial_x^2 f)^2 \big]
	&=\psi_{xx} \int_D (\psi_y^2 f + \partial_x^2 f)^2 + 8 \psi_{xx} \int_{D} (\partial_x^2 f)^2\\
	&\geq 8 \psi_{xx} \int_{D} (\partial_x^2 f)^2.
	\end{split}
\end{equation*}
Now, we write
\begin{equation*}
	-2\int_D \psi_x^2 \psi_{xx} (\partial_y f)^2
	=3 \int_{D} \psi_x^2 \psi_{xx} f \partial_y^2 f 
	+ \int_D \psi_x^2 \psi_{xx} (\partial_y f)^2.
\end{equation*}	
The first term in the right-hand side can be seen as part of a squared term, indeed
\begin{equation*}
	\begin{split}
	\psi_{xx} \int_D \big[ 9\psi_x^4 f^2 + 3 \psi_x^2 f \partial_y^2 f + (\partial_y^2 f)^2\big]
	&=\psi_{xx} \int \big(\tfrac{3}{2} \psi_x^2 f + \partial_y^2 f\big)^2 + \frac{27}{4}\psi_{xx}\int_{D} \psi_x^4 f^2\\
	&\geq \frac{27}{4} \psi_{xx}\int_{D} \psi_x^4 f^2.
	\end{split}
\end{equation*}  
These two estimates together give
\begin{equation}\label{eq:Carl-last}
\begin{split}
		\langle [\mathcal{S}, \mathcal{A}]f,f \rangle
		\ge &\int_D (\tfrac{27}{4}\psi_x^4 \psi_{xx}+ 6\psi_{xt}\psi_{x}^2
		+\psi_{tt} + \tfrac{26}{9} \psi_x^2 \psi_y^2 \psi_{xx} + 2 \psi_{xt}\psi_y^2
		)f^2\\
		&+\int_D (2\psi_x^2\psi_{xx}-6\psi_{xt}) (\partial_x f)^2\\
		&+\int_D (\psi_y^2 \psi_{xx} +\psi_x^2 \psi_{xx}-2\psi_{xt})(\partial_y f)^2\\
		&+ \psi_{xx}\int_D (\partial_y\partial_x f)^2
		+8\psi_{xx}\int_D (\partial_x^2 f)^2.
	\end{split}
\end{equation}
We observe now that 
\begin{equation*}
	\begin{split}
	\psi_x=\frac{2\alpha}{R}\left(\frac{x}{R} + \varphi(t) \right);\qquad
	&\psi_y=\frac{2\alpha y}{R^2};\qquad
	\psi_t=2\alpha \left(\frac{x}{R} + \varphi(t) \right) \varphi'(t);\\
	\psi_{tt}=2\alpha \left(\frac{x}{R} + \varphi(t) \right) \varphi''(2) + 2\alpha &(\varphi'(t))^2;\qquad
	\psi_{xt}= \frac{2\alpha}{R}\varphi'(t); \qquad
	\psi_{xx}=\psi_{yy}=\frac{2\alpha}{R^2}.
	\end{split}
\end{equation*}
Using the above information we can write more explicitly~\eqref{eq:Carl-last}, that is:
\begin{equation*}
	\begin{split}
		\langle [\mathcal{S}, \mathcal{A}]f,f \rangle
		{\color{red}\ge}&
		\int_{D}\left[ 
		216 \frac{\alpha^5}{R^6} \left(\frac{x}{R}+ \varphi(t) \right)^4
		+48 \frac{\alpha^3}{R^3}\varphi'(t)	\left(\frac{x}{R}+ \varphi(t) \right)
		+ 2\alpha \left( \frac{x}{R} + \varphi(t)\right)\varphi''(t)	\right.\\
		&\phantom{\int_D[}\left.+2\alpha(\varphi'(t))^2
		+\frac{208}{9}\frac{\alpha^3}{R^4}\left(\frac{x}{R} + \varphi(t) \right)^2\psi_y^2 + \frac{4\alpha}{R}\varphi'(t)\psi_y^2
		\right]f^2\\
		& +\int_D \left[ 16\frac{\alpha^3}{R^4}\left(\frac{x}{R}+ \varphi(t) \right)^2 -12 \frac{\alpha}{R}\varphi'(t) \right](\partial_x f)^2\\
		&+\int_D \left[\frac{2\alpha}{R^2}\psi_y^2+ \frac{8\alpha^3}{R^4}\left(\frac{x}{R} + \varphi(t) \right)^2 -\frac{4\alpha}{R}\varphi'(t) \right](\partial_y f)^2\\
		&+\frac{2\alpha}{R^2}\int_D (\partial_y\partial_x f)^2 
		+ 16\frac{\alpha}{R^2}\int_D (\partial_x^2 f)^2.
	\end{split}
\end{equation*}
Since $\alpha\geq \overline{C}R^{3/2}>1,$ with $\overline{C}:=\max\{\|\varphi'\|_{L^\infty}, \|\varphi''\|_{L^\infty},1\}$ and since $|\tfrac{x}{R}+ \varphi(t)|\geq 1$ one easily has
\begin{multline*}
	216 \frac{\alpha^5}{R^6} \left(\frac{x}{R}+ \varphi(t) \right)^4
		+48 \frac{\alpha^3}{R^3}\varphi'(t)	\left(\frac{x}{R}+ \varphi(t) \right)
		+ 2\alpha \left( \frac{x}{R} + \varphi(t)\right)\varphi''(t)\\
		\begin{aligned}
		&\geq 216 \frac{\alpha^5}{R^6}\left(\frac{x}{R}+ \varphi(t) \right)^4
		-48\frac{\alpha^3}{R^3}\|\varphi'\|_{L^\infty}\left|\frac{x}{R}+ \varphi(t) \right|
		-2\alpha \|\varphi''\|_{L^\infty}\left|\frac{x}{R}+ \varphi(t) \right|\\
		&\geq 166 \frac{\alpha^5}{R^6}\left(\frac{x}{R}+ \varphi(t) \right)^4. 	
		\end{aligned}
\end{multline*}
Moreover
\begin{equation*}
	\begin{split}
		\frac{208}{9}\frac{\alpha^3}{R^4}\left(\frac{x}{R} + \varphi(t) \right)^2\psi_y^2 + \frac{4\alpha}{R}\varphi'(t)\psi_y^2
		&\geq \frac{208}{9}\frac{\alpha^3}{R^4}\left(\frac{x}{R} + \varphi(t) \right)^2\psi_y^2 
		- \frac{4\alpha}{R}\|\varphi'\|_{L^\infty}\psi_y^2\\
		&\geq \frac{172}{9}\frac{\alpha^3}{R^4}\left(\frac{x}{R} + \varphi(t) \right)^2 \psi_y^2.
	\end{split}
\end{equation*}
Similarly
\begin{equation*}
	\begin{split}
		16\frac{\alpha^3}{R^4}\left(\frac{x}{R}+ \varphi(t) \right)^2 -12 \frac{\alpha}{R}\varphi'(t)
		&\geq 16\frac{\alpha^3}{R^4}\left(\frac{x}{R}+ \varphi(t) \right)^2 
		-12 \frac{\alpha}{R}\|\varphi'\|_{L^\infty}\\
		&\geq \frac{4\alpha^3}{R^4} \left(\frac{x}{R} + \varphi(t)\right)^2\\
		&\geq \frac{4\alpha^3}{R^4},
	\end{split}
\end{equation*}
and 
\begin{equation*}
	\begin{split}
		\frac{8\alpha^3}{R^4}\left(\frac{x}{R} + \varphi(t) \right)^2 -\frac{4\alpha}{R}\varphi'(t)
		&\geq \frac{8\alpha^3}{R^4}\left(\frac{x}{R} + \varphi(t) \right)^2 
		-\frac{4\alpha}{R}\|\varphi'\|_{L^\infty}\\
		&\geq \frac{4\alpha^3}{R^4} \left(\frac{x}{R} + \varphi(t)\right)^2\\
		&\geq \frac{4\alpha^3}{R^4}.
	\end{split}
\end{equation*}
From these estimates and discarding some positive term we get
\begin{equation*}
	\begin{split}
		\langle [\mathcal{S}, \mathcal{A}]f, f\rangle
		\geq & \int_D \left[ 166 \frac{\alpha^5}{R^6}\left(\frac{x}{R}+ \varphi(t) \right)^4 + \frac{172}{9}\frac{\alpha^3}{R^4}\left(\frac{x}{R} + \varphi(t) \right)^2 \psi_y^2  \right] f^2\\
		&+\int_D \frac{4\alpha^3}{R^4} \big[(\partial_x f)^2 + (\partial_y f)^2\big].
	\end{split}
\end{equation*}
Noticing that
\begin{equation*}
	166\frac{\alpha^5}{R^6} \left( \frac{x}{R} + \varphi(t)\right)^4
	\geq134\frac{\alpha^5}{R^6} 
	+ \frac{4\alpha^3}{R^4} (\psi_{xx} + \psi_{yy}) 
	+ \frac{4\alpha^3}{R^4}\psi_x^2
\end{equation*}
and that 
\begin{equation*}
	\frac{172}{9}\frac{\alpha^3}{R^4} \left(\frac{x}{R} + \varphi(t) \right)^2\psi_y^2\geq \frac{4\alpha^3}{R^4} \psi_y^2,
\end{equation*}
one gets
\begin{equation}\label{eq:last}
	\langle [\mathcal{S}, \mathcal{A}]f, f \rangle
	\geq 134\frac{\alpha^5}{R^6} \int_D f^2 
	+\frac{4\alpha^3}{R^4}\int_{D} \big[(\partial_x f)^2 + (\psi_{xx} + \psi_x^2)f^2\big] 
	+\frac{4\alpha^3}{R^4}\int_{D} \big[(\partial_y f)^2 + (\psi_{yy} + \psi_y^2)f^2\big]. 
\end{equation}
Observing that 
\begin{equation*}
	\int_{D} [(\partial_x -\psi_x)f]^2
	=\int_D \big[(\partial_x f)^2 + (\psi_{xx} + \psi_x^2)f^2 \big],
\end{equation*}
and similarly for the corresponding term in the $y$ variable, from~\eqref{eq:last} one obtains~\eqref{eq:Carl-conj}.
\end{proof}

\section*{Acknowledgements}

L.C was supported by the grant Ramón y Cajal RYC2021-032803-I  funded by MCIN/AEI/10.13039/50110
\noindent
0011033 and by the European Union NextGenerationEU/PRTR, the Deutsche Forschungsgemeinschaft (DFG, German
Research Foundation) -- Project-ID 258734477 -- SFB 1173 and by Ikerbasque.

L.F. was supported by project PID2021-123034NB-I00/AEI/10.13039/501100011033 funded by Agencia de Investigaci\'on, by the project IT1615-22 funded by the Basque Government, and by Ikerbasque.

F.L. was partially supported by CNPq grant 305791/2018-4 and FAPERJ grant E-26/202.638/2019.

\providecommand{\bysame}{\leavevmode\hbox to3em{\hrulefill}\thinspace}


\end{document}